\documentclass[12pt]{article}

\usepackage{amssymb,amsmath,latexsym,graphics,supertabular}
\usepackage[enableskew]{youngtab}
\usepackage{amsbsy}
\usepackage{amscd}
\usepackage{verbatim}
\usepackage{amsthm}
\usepackage{mathrsfs}
\usepackage{verbatim}
\usepackage{graphicx,subfigure}
\usepackage[colorlinks]{hyperref}
\usepackage{xypic}
\usepackage{url}
\usepackage{tikz}
    \usetikzlibrary{backgrounds,fit,decorations.pathreplacing}
    \usetikzlibrary{arrows,shapes,positioning}
    \usetikzlibrary{calc,through,backgrounds,chains}
    \usetikzlibrary{arrows,shapes,snakes,automata,backgrounds,petri}
    \usetikzlibrary{matrix,arrows,decorations.pathmorphing}
\usepackage{todonotes}
\usepackage{authblk}

\newenvironment{pf*}[1]{\proof[#1]}{\endproof}
\newtheorem{Theorem}[equation]{Theorem}
\newtheorem*{Theorem2}{Theorem}
\newtheorem{Corollary}[equation]{Corollary}
\newtheorem{Lemma}[equation]{Lemma}
\newtheorem{Proposition}[equation]{Proposition}

\theoremstyle{definition}
\newtheorem{Definition}[equation]{Definition}
\newtheorem{Algorithm}[equation]{Algorithm}
\newtheorem{Example}[equation]{Example}

\theoremstyle{remark}

\newtheorem{Exercise}[equation]{Exercise}
\newtheorem{Remark}[equation]{Remark}

\numberwithin{equation}{section}
\numberwithin{figure}{section}

\newcommand{\beqn}{\begin{equation}}
\newcommand{\eeqn}{\end{equation}} \newcommand{\bdfn}{\begin{Definition}}
\newcommand{\edfn}{\end{Definition}}
\newcommand{\brem}{\begin{Remark}}
\newcommand{\erem}{\end{Remark}}
\newcommand{\benum}{\begin{enumerate}}
\newcommand{\eenum}{\end{enumerate}}
\newcommand{\bexam}{\begin{Example}}
\newcommand{\eexam}{\end{Example}}
\newcommand{\bexer}{\begin{Exercise}}
\newcommand{\eexer}{\end{Exercise}}
\newcommand{\bthm}{\begin{Theorem}}
\newcommand{\ethm}{\end{Theorem}}

\newcommand{\blem}{\begin{Lemma}}
\newcommand{\elem}{\end{Lemma}}
\newcommand{\bprop}{\begin{Proposition}}
\newcommand{\eprop}{\end{Proposition}}
\newcommand{\bcor}{\begin{Corollary}}
\newcommand{\ecor}{\end{Corollary}}
\newcommand{\beqna}{\begin{equation*}}
\newcommand{\eeqna}{\end{equation*}}

\begin{document}

\title{Chains in Weak Order Posets Associated to Involutions}

\author[1]{Mahir Bilen Can}
\author[2]{Michael Joyce}
\author[3]{Benjamin Wyser}
\affil[1]{{\small Tulane University, New Orleans; mcan@tulane.edu}}
\affil[2]{{\small Tulane University, New Orleans; mjoyce3@tulane.edu}}
\affil[3]{{\small University of Illinois at Urbana-Champaign; bwyser@illinois.edu}}
\normalsize

\maketitle

\begin{abstract}
The $\mathcal{W}$-set of an element of a weak order poset is useful in the cohomological study of the closures of spherical subgroups in generalized flag varieties.  We explicitly describe in a purely combinatorial manner the $\mathcal{W}$-sets of the weak order posets of three different sets of involutions in the symmetric group, namely, the set of all involutions, the set of all fixed point free involutions, and the set of all involutions with signed fixed points (or ``clans").  These distinguished sets of involutions parameterize Borel orbits in the classical symmetric spaces associated to the general linear group.  In particular, we give a complete characterization of the maximal chains of an arbitrary lower order ideal in any of these three posets.
\end{abstract}

\section{Introduction}\label{sec: intro}

Given a reductive group $G$, an algebraic variety $X$ equipped with a $G$-action is said to be spherical if a Borel subgroup $B$ of $G$ has a dense orbit in $X$.  (All varieties in this paper are defined over an algebraically closed field of characteristic $\neq 2$.)  A subgroup $H\subseteq G$ is called spherical if the homogeneous space $G / H$ is spherical.  The geometry of spherical varieties provides a rich source of combinatorial structures.

A spherical $G$-variety $X$ always has finitely many $B$-orbits \cite{Brion86, Popov86, VK78}.
We denote the set of $B$-orbit closures in $X$ by $\mathcal{B}(X)$. The set $\mathcal{B}(X)$ possesses two geometrically natural partial orders:  the Bruhat order (given by inclusion), and the weak order.  The {\em weak order} is the transitive closure of the covering relations given by $Y_1 \lessdot Y_2$ if and only if $Y_2 = P Y_1$ for some minimal parabolic subgroup $P$ containing $B$. Not much is known about these poset structures for arbitrary spherical varieties -- a parametrization of the $B$-orbits is unknown in general.  Yet, in the special case of symmetric homogeneous spaces, there is a natural order-preserving map from the weak order poset $\mathcal{B}(X)$ to the weak order poset of twisted involutions of the Weyl group $W$ of $G$ \cite{RS90}.  Parameterizations of $\mathcal{B}(X)$ for all classical type symmetric homogeneous spaces are given in \cite{MO90}.  See \cite{RS93} and \cite {Helminck04} for more on the combinatorics of Borel orbits of symmetric varieties.

In this paper, we study maximal chains in the weak order on three sets of involutions:
\begin{enumerate}
 \item The poset of all involutions in the symmetric group $S_n$.
 \item The poset of all fixed-point free involutions in $S_n$.
 \item The poset of all involutions with signed fixed points of constant total charge.  The elements of this poset are often referred to as $(p,q)$-clans, where $p+q=n$.  The precise definition is given in Section \ref{sec: sgn inv W-set}.
\end{enumerate}
These posets are the opposite of weak order posets $\mathcal{B}(X)$ for three classical symmetric homogeneous spaces $G/H$ of type $A$, where $G=\text{GL}_n$, and $H$ is, respectively, the central extension of the orthogonal group $\text{O}_n$, the central extension of the symplectic group $\text{Sp}_n$, and the product subgroup $\text{GL}_p \times \text{GL}_q$.  We combinatorially characterize chains in the three weak order posets of involutions listed above.  Such descriptions are useful for understanding the stratification of $H$-orbits in $G/B$, or equivalently of $B$-orbits in $G/H$ and constitute one of the fundamental problems in the study of spherical varieties.  While our results have geometric significance, we emphasize that the results and methods in this paper are purely combinatorial in nature and do not rely on any geometric considerations.  The results obtained here will be combined with geometric arguments to obtain new Schubert polynomial identities in a future work \cite{CJW14b}.

Our primary combinatorial object of study is the $\mathcal{W}$-set of an element of one of our three posets. The notion of $\mathcal{W}$-set is introduced by Brion in \cite{Brion01} in a geometric context.  Here, we give a purely combinatorial definition.  Let $P$ be a poset with each covering relation $\lessdot$ assigned a set of labels from the set $\{1, 2, \dots, n-1\}$.  If $p \lessdot p'$ and $j$ belongs to the set of labels for $\lessdot$, write $p \lessdot_j p'$.  Assuming $P$ to have a unique minimal element $\hat{0}$, the \textbf{$\mathcal{W}$-set} of an element $p \in P$, denoted $\mathcal{W}(p)$, consists of all $w = s_{j_\ell} \cdots s_{j_2} s_{j_1} \in S_n$ of length $\ell$ such that
\[ \hat{0} = p_0 \leq_{j_1} p_1 \leq_{j_2} \cdots \leq_{j_\ell} p_\ell = p, \]
for some $p_1, p_2, \dots, p_{\ell-1} \in P$.  (As usual, $s_j$ denotes the simple transposition of $S_n$ which interchanges $j$ and $j+1$.)  In our examples and others arising from the same sort of geometric construction, it follows from results of Richardson and Springer \cite{RS90, RS93} that the maximal chains in any lower order ideal  $[\hat{0}, p] \subseteq P$ are then parameterized by the reduced expressions of the elements of $\mathcal{W}(p)$.

Our main theorems give explicit combinatorial descriptions of the $\mathcal{W}$-sets of arbitrary elements in each of our three involution posets.
Despite their geometric and combinatorial significance, explicit descriptions of $\mathcal{W}$-sets are known in only a small number of cases.  First, if $X$ is any spherical variety, with $Y \in \mathcal{B}(X)$, one can define a rank associated to $Y$, and then speak of whether $Y$ is of ``maximal rank''.  In \cite{Brion01}, Brion gives a parametrization of the $B$-orbit closures $Y$ on $X$ of maximal rank, as well as a description of $\mathcal{W}(Y)$ for any such $Y$.  For a spherical variety in which all $B$-orbit closures have maximal rank (such a space is said to be of ``minimal rank''), this gives complete information.  The spherical homogeneous spaces of minimal rank are determined by Ressayre in \cite{Ressayre10}.

Aside from cases covered by the above, the only explicit determination of a $\mathcal{W}$-set in the literature of which we are aware is given in \cite{CJ13}.  The primary result of that paper describes the $\mathcal{W}$-set of the longest permutation in the poset of all involutions of $S_n$.  Thus the results of this paper are a considerable generalization of the main result of \cite{CJ13}.

The organization of the paper is as follows.  After giving some preliminary notation and definitions, we describe the $\mathcal{W}$-set of any element of the weak order poset for involutions in Section \ref{sec: inv W-set}.  Then in Section \ref{sec: matchings}, a proof of this description is given using the language of matchings.  In Section \ref{sec: fpf inv W-set}, we deduce a description of the $\mathcal{W}$-set of any element of the weak order poset for fixed-point free involutions from the result for involutions.  In Section \ref{sec: sgn inv W-set}, we describe weak order posets of clans, and characterize the $\mathcal{W}$-sets of their elements.

\section{Involution combinatorics}\label{sec: inv comb}

\subsection{Preliminaries}\label{sec: prelim}
As usual, $[n]$ denotes the set $\{1,\hdots,n\}$.  We write the elements of the symmetric group $S_n$ in cycle notation using parentheses, as well as in one-line notation using brackets, typically reserving cycle notation for permutations that belong to one of our involution posets.  We omit brackets in one-line notation when no ambiguity is present.  For example, $w=4213=[4,2,1,3]=(1,4,3)$ is the permutation sending $1$ to $4$, $2$ to $2$, $3$ to $1$, and $4$ to $3$. Accordingly, we multiply elements of $S_n$ by applying composition from right to left. The identity permutation of a symmetric group is denoted by $\text{id}$.

If $w \in S_n$ and $i, j \in [n]$, we say that {\em $i$ occurs before $j$ in $w$} if $w^{-1}(i) < w^{-1}(j)$, or in other words, if $i$ occurs to the left of $j$ in the one-line notation for $w$.  We similarly speak of {\em $i$ occurring after $j$} and {\em $k$ occurring between $i$ and $j$} in $w$.  If $A \subseteq [n]$ and $i,j \in A$, then $i$ and $j$ are said to be {\em adjacent in $A$} if there is no $k \in A$ that occurs between $i$ and $j$ in $A$.

Recall that the symmetric group $S_n$ is a Coxeter group with simple generators $S := \{s_1,s_2,\dots,s_{n-1}\}$, where $s_i = (i, i+1)$ is the simple transposition that interchanges $i$ and $i+1$.  The relations that define $S_n$ are
$$
s_i^2 = \text{id} \text{ for all $1 \leq i \leq n-1$}
$$
and the braid relations,
\begin{equation}\label{eqn:braid rels}
\begin{cases}
s_i s_j = s_j s_i & \text{if $|j - i| > 1$,} \\
s_i s_{i+1} s_i = s_{i+1} s_i s_{i+1} & \text{for all $1 \leq i \leq n-1$.}
\end{cases}
\end{equation}
If $w = w_1 w_2 \cdots w_k$ is a word from an ordered alphabet, then an inversion of $w$ is a pair $i <j$ such that $w_i > w_j$.  If $w \in S_n$, then the length $\ell(w)$ of $w$ is the number of inversions of the word corresponding to the one-line notation of $w$.  It is also equal to the least value of $k$ such that $w = s_{i_1} \cdots s_{i_k}$ for some $s_{i_1}, \dots, s_{i_k} \in S$.

The \emph{Richardson-Springer monoid} associated to $S_n$, which we denote by $M(S_n)$, is the finite monoid generated
by $S' := \{m(s_1), m(s_2), \dots, m(s_{n-1})\}$, subject to the relations $m(s_i)^2 = m(s_i)$ for all
$1 \leq i \leq n-1$ and the braid relations \eqref{eqn:braid rels} with $m(s_i)$ in place of $s_i$.  If for any
$w \in S_n$ we define $m(w) := m(s_{i_1}) \cdots m(s_{i_l})$ for any reduced expression $s_{i_1} \cdots s_{i_l}$
of $w$, then the correspondence $w \leftrightarrow m(w)$ is a well-defined set-theoretic bijection
$S_n \leftrightarrow M(S_n)$ \cite{RS90}.

\subsection{Involutions}\label{sec: inv W-set}

We proceed with defining the weak order on the set of involutions of $S_n$ and then describe the maximal chains of certain intervals therein.

Let $\mathscr{I}_n = \{ \pi \in S_n : \pi^2 = \text{id} \}$ denote the set of all involutions of $S_n$.  It is an $M(S_n)$-set through the action defined inductively by
\begin{equation}\label{eqn:RS mon action on I_n}
m(s_i) \cdot \pi =
\begin{cases}
s_i \pi s_i^{-1} & \text{if $\ell(s_i \pi s_i^{-1}) = \ell(\pi) + 2$,} \\
s_i \pi & \text{if $\pi(i) = i$ and $\pi(i+1) = i+1$,} \\
\pi & \text{otherwise.}
\end{cases}
\end{equation}
We write $\pi \lessdot_{i} \pi'$ if $\pi' = m(s_i) \cdot \pi$ and $\pi' \neq \pi$.  Moreover, we refer to the case where $\pi \lessdot_{i} \pi'$ and $\ell(\pi') = \ell(\pi) + 2$ as a {\em covering relation of type I} and the case where $\pi \lessdot_{i} \pi'$ and $\ell(\pi') = \ell(\pi) + 1$ as a {\em covering relation of type II}.  When we interpret the poset of involutions in terms of matchings, the covering relations of type I will be further broken down into subcases; see Figure \ref{fig:matching types}.

\bdfn
Let $\pi, \pi' \in \mathscr{I}_n$.  The \emph{weak order} on $\mathscr{I}_n$ is the partial order defined by $\pi \leq \pi'$ if and only if there exist $\pi_0 = \pi, \pi_1, \pi_2, \dots, \pi_k = \pi'$ such that for each $1 \leq j \leq k$, $\pi_{j-1} \lessdot_{i} \pi_{j}$ for some $i$.  Equivalently, $\pi' = m(w) \cdot \pi$ for some (not necessarily unique) $m(w) \in M(S_n)$.
\edfn

The Hasse diagram for weak order for $\mathscr{I}_4$ is shown in Figure \ref{fig:weak order I_4}.  Single edges represent pairs $(\pi, \pi')$ such that $\pi' = m(s_i) \cdot \pi = s_i \pi s_i^{-1}$ (type I), while double edges represent pairs $(\pi, \pi')$ such that $\pi' = m(s_i) \cdot \pi = s_i \pi$ (type II).  The geometric significance of this distinction is explained in \cite[Section 1]{Brion01}, but the distinction between single and double edges will not play a role in the sequel.  The poset $\mathscr{I}_n$ has a minimal element $\alpha_n := \text{id}$, and a maximal element $\beta_n := w_0$, where $w_0 = [n, n-1, \cdots, 2, 1]$, the longest permutation.  We drop the subscript $n$ when it is clear from context.

\begin{figure}[htp]
\begin{center}

\begin{tikzpicture}[scale=.4]

\node at (0,0) (a) {$\text{id}$};

\node at (-10,5) (b1) {$(12)$};
\node at (0,5) (b2) {$(23)$};
\node at (10,5) (b3) {$(34)$};

\node at (-10,10) (c1) {$(13)$};
\node at (0,10) (c2) {$(12)(34)$};
\node at (10,10) (c3) {$(24)$};

\node at (-5,15) (d1) {$(14)$};
\node at (5,15) (d2) {$(13)(24)$};

\node at (0,20) (e) {$(14)(23)$};

\node at (-6.5,2.5) {$1$};
\node at (0.5,2.5) {$2$};
\node at (6.5,2.5) {$3$};

\node at (-10.5,7.5) {$2$};
\node at (-7,5.8) {$3$};
\node at (-3,5.8) {$1$};
\node at (3,5.8) {$3$};
\node at (7,5.8) {$1$};
\node at (10.5,7.5) {$2$};

\node at (-8,12.5) {$3$};
\node at (1,11.8) {$2$};
\node at (6,10.8) {$1$};

\node at (-3.5,17.5) {$2$};
\node at (4,17.5) {$1,3$};

\draw[-, double, thick] (a) to (b1);
\draw[-, double, thick] (a) to (b2);
\draw[-, double, thick] (a) to (b3);

\draw[-, thick] (b1) to (c1);
\draw[-, double, thick] (b1) to (c2);
\draw[-, thick] (b2) to (c1);
\draw[-, thick] (b2) to (c3);
\draw[-, double, thick] (b3) to (c2);
\draw[-, thick] (b3) to (c3);

\draw[-, thick] (c1) to (d1);
\draw[-, thick] (c2) to (d2);
\draw[-, thick] (c3) to (d1);

\draw[-, double, thick] (d1) to (e);
\draw[-, thick] (d2) to (e);

\end{tikzpicture}

\caption{Weak order on $\mathscr{I}_4$.}\label{fig:weak order I_4}

\end{center}
\end{figure}

An involution $\pi = (a_1, b_1)(a_2, b_2) \cdots (a_k, b_k)$, expressed as a product of disjoint transpositions, is in \emph{standard form} if $a_i < b_i$ for all $1 \leq i \leq k$ and $a_1 < a_2 < \dots < a_k$.
For $\pi = (a_1, b_1) \cdots (a_k, b_k) \in \mathscr{I}_n$, define $L(\pi)$ by
$$
L(\pi) := \frac{\ell(\pi)+k}{2},
$$
where $\ell(\pi)$ is the length of $\pi$ as an element of $S_n$ and $k$ is the number of disjoint $2$-cycles that appear in the cycle decomposition of $\pi$.  Weak order on $\mathscr{I}_n$ is a ranked poset, with rank function $L$ \cite{RS90}.

\bdfn
The $\mathcal{W}$-set of $\pi \in \mathscr{I}_n$ is
$$
\mathcal{W}(\pi) := \{ w \in S_n : m(w) \cdot \alpha = \pi \text{ and } \ell(w) = L(\pi) \}.
$$
\edfn

It follows immediately from the definition of $M(S_n)$ that this definition of $\mathcal{W}(\pi)$ is the same as that given in introduction when the poset in question possesses a $M(S_n)$-action.  The maximal chains in the interval $[\alpha, \pi]$ are parameterized by the reduced decompositions of the elements of $\mathcal{W}(\pi)$.

We now state our first main combinatorial result.  It is proven in the next subsection where it is reformulated in terms of matchings.

\bthm\label{thm:W-set for I_n}
Let $\pi = (a_1, b_1)(a_2, b_2) \cdots (a_k, b_k) \in \mathscr{I}_n$ be written in standard form, with fixed points $c_1 < c_2 < \dots < c_l$.
In this case, $\mathcal{W}(\pi)$ consists of all $w = [w(1), w(2), \dots, w(n)]$ such that
\benum
\item for each $1 \leq i \leq k$, $b_i$ occurs before $a_i$ in $w$ and for any $a_i < x < b_i$, $x$ does not occur between $b_i$ and $a_i$ in $w$;\label{eq:inv cond 1}
\item if $i < j$ and $b_i < b_j$, then $a_i$ occurs before $b_j$ in $w$;\label{eq:inv cond 2}
\item if $i < j$, then $c_i$ occurs before $c_j$ in $w$;\label{eq:inv cond 3}
\item if $c_j < a_i$, then $c_j$ occurs before $b_i$ in $w$;\label{eq:inv cond 4}
\item if $b_i < c_j$, then $a_i$ occurs before $c_j$ in $w$.\label{eq:inv cond 5}
\eenum
\ethm

\brem
Condition \eqref{eq:inv cond 1} says that for each pair $a_i < b_i$ of numbers that form a 2-cycle of $\pi$, the number $b_i$ occurs to the left of $a_i$ in the one-line notation of $w$.  Moreover, conditions \eqref{eq:inv cond 4} and \eqref{eq:inv cond 5} in conjunction with condition \eqref{eq:inv cond 1} tell us that the value of any fixed point $c_j$ cannot occur between $b_i$ and $a_i$.  Condition \eqref{eq:inv cond 2} tells us the relative order of the pairs ``$b_i \hdots a_i$'' and ``$b_j \hdots a_j$'' unless the values are ordered $a_i < a_j < b_j < b_i$.  In that latter case, the order of the four values in the one-line notation of $w$ can be any of three possibilities: (1) $b_i \hdots a_i \hdots b_j \hdots a_j$; (2) $b_j \hdots a_j \hdots b_i \hdots a_i$; (3) $b_j \hdots b_i \hdots a_i \hdots a_j$.  These three cases can be seen in the calculation of $\mathcal{W}((1,4)(2,3))$ in Example \ref{ex:inv W-set}.  We refer to the order appearing in case (3) as the pair $b_j \hdots a_j$ ``nesting around" $b_i \hdots a_i$.

In order to enumerate the elements of $\mathcal{W}(\pi)$, one may proceed by first placing the values $b_i$ and $a_i$ in order of increasing $i$.  Thus, one first places $b_1$ and $a_1$; since $\pi$ is in standard form, $a_1 = 1$.  The values $b_1$ and $a_1$ cannot nest around any other pair of values from another 2-cycle, so they must occur immediately next to each other in $w$.  So imagining $n$ placeholders arranged horizontally, with the $i^{\text{th}}$ placeholder set to receive the value of $w(i)$, place $b_1$ and $a_1$ in two adjacent places, leaving $n-2$ available placeholders.  Then consider how to place $b_2$ and $a_2$.  Again, the nesting criterion forces us to place $b_2$ and $a_2$ in adjacent available placeholders.  However, conditions \eqref{eq:inv cond 1} and \eqref{eq:inv cond 2} may restrict where $b_2$ and $a_2$ can be placed.  We continue to place successive values $b_i$ and $a_i$, stopping at any point when our placement violates one of the first two conditions.  Note that there are $(n-1)(n-3) \cdots (n - (2k-1))$ possible placements of the $k$ pairs to consider, although, in practice, many of them can be eliminated from consideration at once early in this recursive process.  Once those $2k$ values are placed, the remaining fixed points have to be placed in increasing order by condition \eqref{eq:inv cond 3}.  Then, having constructed a permutation $w$ which satisfies conditions \eqref{eq:inv cond 1} - \eqref{eq:inv cond 3}, we simply compare the ordering of the fixed points relative to certain pairs to verify conditions \eqref{eq:inv cond 4} and \eqref{eq:inv cond 5}.
\erem

\bexam\label{ex:inv W-set}
In $S_4$, $S_5$, and $S_8$, respectively, we have
\begin{align*}
\mathcal{W}((1,4)(2,3)) & = \{ 3241, 3412, 4132 \}; \\
\mathcal{W}((1,3)(2,5)) & = \{ 31452, 31524 \}; \\
\mathcal{W}((1,6)(3,7)(4,8)) & = \{ 25617384, 26157384, 26173584, 26173845, \\
& \qquad 61257384, 61273584, 61273845 \}.
\end{align*}
\eexam

\subsection{Matchings}\label{sec: matchings}

\begin{figure}[htp]
\begin{center}

\begin{tikzpicture}[scale=.75, transform shape]

\tikzstyle{dot} = [circle, minimum width = 4 pt, fill]

\node[dot] (1) [label=below:$1$] {};
\node[dot] (2) [right=of 1, label=below:$2$] {};
\node[dot] (3) [right=of 2, label=below:$3$] {};
\node[dot] (4) [right=of 3, label=below:$4$] {};
\node[dot] (5) [right=of 4, label=below:$5$] {};

\draw [-] (1.east) -- (2.west);
\draw [-] (2.east) -- (3.west);
\draw [-] (3.east) -- (4.west);
\draw [-] (4.east) -- (5.west);

\draw[-] [in=90, out=90] (1.north) to (3.north);
\draw[-] [in=90, out=90] (2.north) to (5.north);

\end{tikzpicture}

\caption{The matching $\mathcal{M}_{\pi}$ for the involution $\pi = (1,3)(2,5) \in \mathscr{I}_5$.}\label{fig:matching example}

\end{center}
\end{figure}
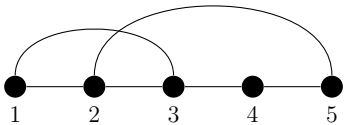

A useful combinatorial model for $\mathscr{I}_n$ is given by the graph theoretic notion of matchings.  A \emph{matching} $\mathcal{M}$ on $n$ vertices consists of a set $\mathcal{V}$ (\emph{isolated vertices}) of singletons and a set $\mathcal{E}$ (\emph{strands}) of doubletons such that $[n]$ is the disjoint union of the singleton sets and the doubleton sets.  The two vertices of a strand are said to be matched.  The matching of an involution $\pi \in \mathscr{I}_n$ is the matching $\mathcal{M}_{\pi}$ whose isolated vertices are the fixed points of $\pi$ and whose strands match distinct vertices $i$ and $j$ if and only if $\pi(i) = j$.  This gives a bijection between $\mathscr{I}_n$, the set of involutions in $S_n$, and the set of all matchings on $n$ vertices.  Figure \ref{fig:matching example} illustrates the matching $\mathcal{M}_{\pi}$ of $\pi = (1,3)(2,5)$.

Two distinct strands $\{ i < j \}$ and $\{ k < l \}$ constitute a \emph{crossing} if either $i < k < j < l$ or $k < i < l < j$.  Two distinct strands $\{ i < j \}$ and $\{ k < l \}$ constitute a \emph{nesting} if either $i < k < l < j$ or $k < i < j < l$.  If $\{ i < j \}$ is a strand, its {\em length} is defined as the difference $j-i$.  The {\em length} $L(\mathcal{M})$ of a matching $\mathcal{M}$ is the sum of the lengths of the strands minus the total number of crossings.

Via the bijection between involutions and matchings, the set of all matchings on $n$ vertices becomes an $M(S_n)$-set.  Write $m(s_i) \cdot \mathcal{M} = \mathcal{M}'$ and $\mathcal{M} \lessdot_{i} \mathcal{M}'$ if $\mathcal{M} = \mathcal{M}_{\pi}$, $\mathcal{M}' = \mathcal{M}_{\pi'}$, and $\pi \lessdot_{i} \pi'$.  Types I and II for covering relations and the weak order carry over to matchings in the obvious manner.  For type I, a cover $\mathcal{M} \lessdot_{i} \mathcal{M}'$ occurs when at most one of the vertices $i$ and $i+1$ is an isolated vertex of $\mathcal{M}$, $\mathcal{M}'$ is obtained from $\mathcal{M}$ by interchanging the strands at vertices $i$ and $i+1$, and the length of $\mathcal{M}'$ is greater than that of $\mathcal{M}$.  See Figure \ref{fig:matching types} for a schematic illustration of all of the possible cases of the action that can occur in type I.  For type II, a cover $\mathcal{M} \lessdot_{i} \mathcal{M}'$ occurs when the vertices $i$ and $i+1$ are both isolated vertices of $\mathcal{M}$, and $\mathcal{M}'$ is obtained from $\mathcal{M}$ by adding a strand that matches vertices $i$ and $i+1$.  This type is illustrated as the final type in Figure \ref{fig:matching types}.

\begin{figure}[htp]
\begin{center}

\begin{tikzpicture}[scale=.75, transform shape]

\tikzstyle{dot} = [circle, minimum width = 4 pt, fill]

\node (1) {};
\node[dot] (2) [right=of 1, label=below:$i$] {};
\node (3) [right=of 2] {};
\node[dot] (4) [right=of 3, label=below:$i+1$] {};
\node (5) [right=of 4] {};

\node (u1) [above=of 1] {};

\draw [-] (1.east) -- (2.west);
\draw [-] (2.east) -- (4.west);
\draw [-] (4.east) -- (5.west);

\draw[-] [in=0, out=90] (2.north) to (u1.east);

\node (6) [right=of 5] {};
\node (7) [right=of 6] {};
\node (8) [right=of 7] {};

\draw[->] (6.east) -- (8.west) node [draw=none,midway,above=1mm] {Type IA1};

\node (9) [right=of 8] {};
\node[dot] (10) [right=of 9, label=below:$i$] {};
\node (11) [right=of 10] {};
\node[dot] (12) [right=of 11, label=below:$i+1$] {};
\node (13) [right=of 12] {};

\node (u9) [above=of 9] {};

\draw [-] (9.east) -- (10.west);
\draw [-] (10.east) -- (12.west);
\draw [-] (12.east) -- (13.west);

\draw[-] [in=45, out=90] (12.north) to (u9.east);

\end{tikzpicture}

\begin{tikzpicture}[scale=.75, transform shape]

\tikzstyle{dot} = [circle, minimum width = 4 pt, fill]

\node (1) {};
\node[dot] (2) [right=of 1, label=below:$i$] {};
\node (3) [right=of 2] {};
\node[dot] (4) [right=of 3, label=below:$i+1$] {};
\node (5) [right=of 4] {};

\node (u5) [above=of 5] {};

\draw [-] (1.east) -- (2.west);
\draw [-] (2.east) -- (4.west);
\draw [-] (4.east) -- (5.west);

\draw[-] [in=180, out=90] (4.north) to (u5.west);

\node (6) [right=of 5] {};
\node (7) [right=of 6] {};
\node (8) [right=of 7] {};

\draw[->] (6.east) -- (8.west) node [draw=none,midway,above=1mm] {Type IA2};

\node (9) [right=of 8] {};
\node[dot] (10) [right=of 9, label=below:$i$] {};
\node (11) [right=of 10] {};
\node[dot] (12) [right=of 11, label=below:$i+1$] {};
\node (13) [right=of 12] {};

\node (u13) [above=of 13] {};

\draw [-] (9.east) -- (10.west);
\draw [-] (10.east) -- (12.west);
\draw [-] (12.east) -- (13.west);

\draw[-] [in=135, out=90] (10.north) to (u13.west);

\end{tikzpicture}

\begin{tikzpicture}[scale=.75, transform shape]

\tikzstyle{dot} = [circle, minimum width = 4 pt, fill]

\node (1) {};
\node[dot] (2) [right=of 1, label=below:$i$] {};
\node (3) [right=of 2] {};
\node[dot] (4) [right=of 3, label=below:$i+1$] {};
\node (5) [right=of 4] {};

\node (u1) [above=of 1] {};
\node (u5) [above=of 5] {};

\draw [-] (1.east) -- (2.west);
\draw [-] (2.east) -- (4.west);
\draw [-] (4.east) -- (5.west);

\draw[-] [in=0, out=90] (2.north) to (u1.east);
\draw[-] [in=180, out=90] (4.north) to (u5.west);

\node (6) [right=of 5] {};
\node (7) [right=of 6] {};
\node (8) [right=of 7] {};

\draw[->] (6.east) -- (8.west) node [draw=none,midway,above=1mm] {Type IB};

\node (9) [right=of 8] {};
\node[dot] (10) [right=of 9, label=below:$i$] {};
\node (11) [right=of 10] {};
\node[dot] (12) [right=of 11, label=below:$i+1$] {};
\node (13) [right=of 12] {};

\node (u9) [above=of 9] {};
\node (u13) [above=of 13] {};

\draw [-] (9.east) -- (10.west);
\draw [-] (10.east) -- (12.west);
\draw [-] (12.east) -- (13.west);

\draw[-] [in=135, out=90] (10.north) to (u13.west);
\draw[-] [in=45, out=90] (12.north) to (u9.east);

\end{tikzpicture}

\begin{tikzpicture}[scale=.75, transform shape]

\tikzstyle{dot} = [circle, minimum width = 4 pt, fill]

\node (1) {};
\node[dot] (2) [right=of 1, label=below:$i$] {};
\node (3) [right=of 2] {};
\node[dot] (4) [right=of 3, label=below:$i+1$] {};
\node (5) [right=of 4] {};

\node (u5) [above=of 5] {};
\node (uu5) [above=of u5] {};

\draw [-] (1.east) -- (2.west);
\draw [-] (2.east) -- (4.west);
\draw [-] (4.east) -- (5.west);

\draw[-] [in=135, out=90] (2.north) to (u5.west);
\draw[-] [in=180, out=90] (4.north) to (uu5.west);

\node (6) [right=of 5] {};
\node (7) [right=of 6] {};
\node (8) [right=of 7] {};

\draw[->] (6.east) -- (8.west) node [draw=none,midway,above=1mm] {Type IC1};

\node (9) [right=of 8] {};
\node[dot] (10) [right=of 9, label=below:$i$] {};
\node (11) [right=of 10] {};
\node[dot] (12) [right=of 11, label=below:$i+1$] {};
\node (13) [right=of 12] {};

\node (u13) [above=of 13] {};
\node (uu13) [above=of u13] {};

\draw [-] (9.east) -- (10.west);
\draw [-] (10.east) -- (12.west);
\draw [-] (12.east) -- (13.west);

\draw[-] [in=165, out=90] (10.north) to (uu13.west);
\draw[-] [in=150, out=90] (12.north) to (u13.west);

\end{tikzpicture}

\begin{tikzpicture}[scale=.75, transform shape]

\tikzstyle{dot} = [circle, minimum width = 4 pt, fill]

\node (1) {};
\node[dot] (2) [right=of 1, label=below:$i$] {};
\node (3) [right=of 2] {};
\node[dot] (4) [right=of 3, label=below:$i+1$] {};
\node (5) [right=of 4] {};

\node (u1) [above=of 1] {};
\node (uu1) [above=of u1] {};

\draw [-] (1.east) -- (2.west);
\draw [-] (2.east) -- (4.west);
\draw [-] (4.east) -- (5.west);

\draw[-] [in=0, out=90] (2.north) to (uu1.east);
\draw[-] [in=45, out=90] (4.north) to (u1.east);

\node (6) [right=of 5] {};
\node (7) [right=of 6] {};
\node (8) [right=of 7] {};

\draw[->] (6.east) -- (8.west) node [draw=none,midway,above=1mm] {Type IC2};

\node (9) [right=of 8] {};
\node[dot] (10) [right=of 9, label=below:$i$] {};
\node (11) [right=of 10] {};
\node[dot] (12) [right=of 11, label=below:$i+1$] {};
\node (13) [right=of 12] {};

\node (u9) [above=of 9] {};
\node (uu9) [above=of u9] {};

\draw [-] (9.east) -- (10.west);
\draw [-] (10.east) -- (12.west);
\draw [-] (12.east) -- (13.west);

\draw[-] [in=30, out=90] (10.north) to (u9.east);
\draw[-] [in=15, out=90] (12.north) to (uu9.east);

\end{tikzpicture}

\begin{tikzpicture}[scale=.75, transform shape]

\tikzstyle{dot} = [circle, minimum width = 4 pt, fill]

\node (1) {};
\node[dot] (2) [right=of 1, label=below:$i$] {};
\node (3) [right=of 2] {};
\node[dot] (4) [right=of 3, label=below:$i+1$] {};
\node (5) [right=of 4] {};

\draw [-] (1.east) -- (2.west);
\draw [-] (2.east) -- (4.west);
\draw [-] (4.east) -- (5.west);

\node (6) [right=of 5] {};
\node (7) [right=of 6] {};
\node (8) [right=of 7] {};

\draw[->] (6.east) -- (8.west) node [draw=none,midway,above=1mm] {Type II};

\node (9) [right=of 8] {};
\node[dot] (10) [right=of 9, label=below:$i$] {};
\node (11) [right=of 10] {};
\node[dot] (12) [right=of 11, label=below:$i+1$] {};
\node (13) [right=of 12] {};

\draw [-] (9.east) -- (10.west);
\draw [-] (10.east) -- (12.west);
\draw [-] (12.east) -- (13.west);

\draw[-] [in=90, out=90] (10.north) to (12.north);

\end{tikzpicture}

\caption{The types of covers $\mathcal{M} \lessdot_{i} \mathcal{M}'$ for matchings.}\label{fig:matching types}

\end{center}
\end{figure}

\blem
The length of an involution $\pi \in \mathscr{I}_n$ and the length of the corresponding matching $\mathcal{M}_{\pi}$ are the same, $L(\mathcal{M}_{\pi}) = L(\pi)$.
\elem

\begin{proof}
The proof is by induction on $L(\pi)$.  First, $\pi = \text{id}$ has length $0$, as does $\mathcal{M}_{\pi}$, which
is the matching without any strands.  Let $\pi$ be any other involution.  Then there exists an involution $\pi' \neq \pi$ and a simple reflection $s_i$ such that $m(s_i) \cdot \pi' = \pi$, and so also $m(s_i) \cdot \mathcal{M}_{\pi'} = \mathcal{M}_{\pi}$.  Since $L(\mathcal{M}_{\pi'}) = L(\pi')$ by the induction hypothesis, it suffices to show that $L(\pi) = L(\pi') + 1$ and $L(\mathcal{M}_{\pi}) = L(\mathcal{M}_{\pi'}) + 1$.  For the first equality, note that either $\ell(\pi) = \ell(\pi') + 2$ (in type I) or $\ell(\pi) = \ell(\pi') + 1$ and $\pi$ has one more $2$-cycle than $\pi'$ (in type II).  For the second equality, note that a type I covering relation corresponds to either the length of one of the strands increasing by one (type IA), the length of two strands increasing by one but also adding a crossing (type IB), or the length of one strand increasing by one, the length of a second strand decreasing by one, and a crossing being eliminated (type IC).  In a type II covering relation, no crossings are introduced and no strand lengths change, but a new strand of length one is created.
\end{proof}

We now reformulate Theorem \ref{thm:W-set for I_n} in terms of matchings and give a proof.

\bthm\label{thm:W-set for matchings}
If $\mathcal{M} := \mathcal{M}_{\pi}$ denotes the corresponding matching for the involution $\pi \in \mathcal{I}_n$, then $\mathcal{W}(\pi)$ consists of all $w = [w(1),w(2),\dots,w(n)]$ such that all of the following hold for all choices of the various indices:
\benum
\item If $\{ i < j \}$ is a strand of $\mathcal{M}$, then $j$ occurs before $i$ in $w$ and if $k$ is any vertex between $i$ and $j$, then $k$ either occurs before $j$ or after $i$ in $w$ (not in between $j$ and $i$).\label{eq:match cond 1}
\item If $\{ i < j \}$ and $\{ k < l \}$ are two non-nesting strands of $\mathcal{M}$ with $i < k$ and $j < l$, then $i$ occurs before $l$ in $w$.\label{eq:match cond 2}
\item If $i$ and $j$ are isolated vertices of $\mathcal{M}$ with $i < j$, then $i$ occurs before $j$ in $w$.\label{eq:match cond 3}
\item If $\{ i < j \}$ is a strand and $k$ is an isolated vertex of $\mathcal{M}$ with $k < i$, then $k$ occurs before $j$ in $w$.\label{eq:match cond 4}
\item If $\{ i < j \}$ is a strand and $k$ is an isolated vertex of $\mathcal{M}$ with $j < k$, then $i$ occurs before $k$ in $w$.\label{eq:match cond 5}
\eenum
\ethm

\begin{proof}

We show by induction on $L(\pi)$ that $w \in \mathcal{W}(\pi)$ if and only if $w$ satisfies \eqref{eq:match cond 1} -- \eqref{eq:match cond 5}.  If $\pi = \text{id}$, then $\mathcal{M}$ is the matching with no strands, so all of the conditions are vacuously satisfied except \eqref{eq:match cond 3}, which only holds for $w = \text{id} = [1,2,\dots,n]$, the unique element of $\mathcal{W}(\text{id})$.

Let $\pi \in \mathscr{I}_n$ with $L(\pi) > 0$.  First, we show that every $w \in \mathcal{W}(\pi)$ satisfies \eqref{eq:match cond 1} - \eqref{eq:match cond 5}.  Let $w \in \mathcal{W}(\pi)$.  Then, by definition, there exists $\pi' \in \mathscr{I}_n$, $w' \in \mathcal{W}(\pi')$ and $s_i \in S$ such that $\pi = m(s_i) \cdot \pi'$ with $L(\pi') < L(\pi)$ and $w = s_i w'$ with $\ell(w) > \ell(w')$.  In particular, $i$ occurs before $i+1$ in $w'$, their positions are switched in $w$, with $i+1$ occurring before $i$, and all other values of $w'$ unchanged.  By induction, $w'$ satisfies \eqref{eq:match cond 1} - \eqref{eq:match cond 5} relative to $\mathcal{M}' := \mathcal{M}_{\pi'}$.  Checking that $w$ satisfies the same conditions relative to $\mathcal{M}$ amounts to a case-by-case check.  There are a large number of cases to check.  For the sake of brevity, for each type of covering relation and for each of conditions \eqref{eq:match cond 1} - \eqref{eq:match cond 5}, we give the argument that $w$ satisfies the condition only if it is not immediately implied by the fact that $w'$ satisfies that same condition.  The omitted cases, though numerous, are trivial to check.

Suppose $\pi' \lessdot_{i} \pi$ is of type IA1.  Thus $\mathcal{M}'$ contains a strand $\{j < i\}$ and an isolated vertex at $i+1$, while $\mathcal{M}$ contains the strand $\{j < i+1\}$ and isolated vertex $i$.  To show \eqref{eq:match cond 1} holds for $w$, we must verify that $i+1$ occurs before $j$ and no $j < k < i+1$ occurs between $i+1$ and $j$.  That $i+1$ occurs before $j$ follows from the fact that $i$ occurs before $j$ in $w'$, and $i+1$ occupies the same spot in $w$ as $i$ does in $w'$, with $j$ unchanged.  As for the second part of condition \eqref{eq:match cond 1}, it follows immediately for $k \in (j,i)$ from the corresponding condition for $w'$, while the fact that $i$ does not occur between $j$ and $i+1$ in $w$ follows from the fact that $i+1$ occurs before $i$ in $w$, since $w=s_iw' > w'$.  The reasoning for type IA2 is similar.

Now suppose $\pi' \lessdot_{i} \pi$ is of type IB.  Thus $\mathcal{M}'$ contains strands $\{j < i\}$ and $\{i+1 < k\}$, while $\mathcal{M}$ contains strands $\{j<i+1\}$ and $\{i<k\}$.  To show \eqref{eq:match cond 1} holds for $w$, we must show that $i$ does not occur between $j$ and $i+1$ in $w$, and that $i+1$ does not occur between $i$ and $k$ in $w$.  Equivalently, we must show that $i+1$ does not occur between $i$ and $j$ in $w'$, and that $i$ does not occur between $i+1$ and $k$ in $w'$.  Both of these conditions follow from \eqref{eq:match cond 2} for $w'$, since the two strands in question in $\mathcal{M}'$ are non-nesting, which (together with the fact that $s_iw' > w'$) implies that $w'$ must be of the form $\hdots i \hdots j \hdots k \hdots i+1 \hdots$.

Next, suppose that $\pi' \lessdot_{i} \pi$ is of type IC1.  Thus $\mathcal{M}'$ contains strands $\{i < j\}$ and $\{i+1 < k\}$ with $i < i+1 < j < k$, while $\mathcal{M}$ contains nested strands $\{i<k\}$ and $\{i+1<j\}$.  To establish \eqref{eq:match cond 1} for $w$, we must see that neither $i+1$ nor $j$ occurs between $i$ and $k$ in $w$, or equivalently that neither $i$ nor $j$ occurs between $i+1$ and $k$ in $w'$.  Again, this follows from the fact that $w'$ satisfies \eqref{eq:match cond 2}.  Indeed, applying \eqref{eq:match cond 2} to the non-nested strands $\{i<j\}$ and $\{i+1<k\}$, we see that $w'$ must be of the form $\hdots j \hdots i \hdots k \hdots i+1 \hdots$.  The reasoning for type IC2 is similar.

Lastly, suppose $\pi' \lessdot_{i} \pi$ is of type II.  Thus $\mathcal{M}'$ contains isolated vertices $i$ and $i+1$, while $\mathcal{M}$ contains the strand $\{i<i+1\}$.  The second part of condition \eqref{eq:match cond 1} for $w$ holds vacuously, while the first follows from the fact that $i$ occurs before $i+1$ in $w'$, as follows either from the fact that $s_iw' > w'$, or alternatively from condition \eqref{eq:match cond 3} applied to $w'$.  For condition \eqref{eq:match cond 2}, suppose that there is an arc $\{j < k\}$ in $\mathcal{M}$ with $i+1<j<k$.  We must see that $i$ occurs before $k$ in $w$.  This follows from the fact that $i+1$ occurs before $k$ in $w'$, which can be deduced from \eqref{eq:match cond 4} applied to $\mathcal{M}'$.  One verifies similarly that \eqref{eq:match cond 2} applies also to any strand $\{j < k\}$ with $j < k < i$.

To check condition \eqref{eq:match cond 4} for $w$, suppose that $k < i$ is an isolated vertex.  We must show that $k$ occurs before $i$ in $w$, or equivalently that $k$ occurs before $i+1$ in $w'$.  This follows from condition \eqref{eq:match cond 3} applied to $w'$.  Condition \eqref{eq:match cond 5} for $w$ is verified similarly.

We now prove the converse statement.  Suppose that $\pi \in \mathscr{I}_n$, and suppose that $w$ satisfies the five conditions of the theorem for $\mathcal{M}_{\pi}$.  The goal is to find $\pi' \in \mathscr{I}_n$ such that $\pi' \lessdot_{i} \pi$ for some $i$, and then to observe that $w'=s_iw$ satisfies conditions \eqref{eq:match cond 1} - \eqref{eq:match cond 5} of the theorem for $\mathcal{M}_{\pi'}$.  Then by induction, we will have that $w' \in \mathcal{W}(\pi')$, and hence $w \in \mathcal{W}(\pi)$.  Note that in order for this to occur, we must have $\ell(w) > \ell(w')$, or equivalently, $i+1$ must occur before $i$ in $w$.

Now, if $\{ i < j \}$ is a strand of $\mathcal{M}_{\pi}$, call $i$ the left vertex of the strand, and $j$ the right vertex.  There exists a pair of vertices $i < j$ such that $i$ is the left vertex of a strand, $j$ is the right vertex of a strand, and every vertex in between is an isolated vertex.  Indeed, if $\mathcal{M}_{\pi}$ has $e$ strands, construct a word of length $2e$ consisting of $e$ L's and $e$ R's by labeling the vertices of strands L for left vertices and R for right vertices, reading left to right.  This word starts with L and ends with R, so it must contain `LR' in consecutive character positions at least once.

First, consider the case where there are no isolated vertices between $i$ and $j$, so that $j=i+1$.  If $\{ i < i+1 \}$ is a strand of $\mathcal{M}$, then a type II covering relation is available with $\mathcal{M}' \lessdot_i \mathcal{M}$.  Then $i+1$ occurs before $i$ in $w$ by condition \eqref{eq:match cond 1} for $\mathcal{M}$.  We have that $w' = s_i w$ satisfies condition \eqref{eq:match cond 3} for $\mathcal{M}'$ because $w$ satisfies conditions \eqref{eq:match cond 4} and \eqref{eq:match cond 5} for $\mathcal{M}$.  In addition, $w'$ satisfies conditions \eqref{eq:match cond 4} and \eqref{eq:match cond 5} for $\mathcal{M}'$ due to the fact that $w$ satisfies condition \eqref{eq:match cond 2} for $\mathcal{M}$.  Conditions \eqref{eq:match cond 1} and \eqref{eq:match cond 2} for $\mathcal{M}'$ hold for $w'$ by virtue of $w$ satisfying the corresponding conditions for $\mathcal{M}$.

On the other hand, if $i$ and $i+1$ are vertices of crossing strands $\{ i < l \}$ and $\{ k < i +1 \}$ with $k < i < i+1 < l$, then a type IB covering relation, $\mathcal{M}' \lessdot_i \mathcal{M}$ is available.  In this case, each of conditions \eqref{eq:match cond 1} - \eqref{eq:match cond 5} for $w' = s_i w$ follow from the corresponding condition for $w$.

If there is an isolated vertex between $i$ and $j$, then we claim that either $i+1$ occurs before $i$, or $j$ occurs before $j-1$ in $w$.  Indeed, $j$ occurs before $i$ in $w$, either by condition \eqref{eq:match cond 1} if $\{i<j\}$ is a strand of $\mathcal{M}$, or by conditions \eqref{eq:match cond 1} and \eqref{eq:match cond 2} if $\{k<j\}$ and $\{i<l\}$ with $k<i<j<l$ are crossing strands.  Now, if $i+1=j-1$, then either this value occurs to the left of $i$, in which case we are done, or it occurs to the right of $i$, hence also to the right of $j$, in which case we are also done.  On the other hand, if $i+1<j-1$, then by condition \eqref{eq:match cond 3}, $i+1$ must occur left of $j-1$ in $w$.  So if $i+1$ is left of $i$, then we are done, and otherwise $i+1$ is right of $i$, hence right of $j$, and so $j-1$ is also right of $j$.

Suppose that $i+1$ occurs before $i$ in $w$.  We claim that $w' = s_i w$ satisfies conditions \eqref{eq:match cond 1} - \eqref{eq:match cond 5} for $\mathcal{M}'$, where $\mathcal{M}' \lessdot_{i} \mathcal{M}$ is of type IA2.  Indeed, each condition for $w'$ follows from the corresponding condition for $w$.  If instead $j$ occurs before $j-1$ in $w$, then $w' = s_{j-1}w$ satisfies conditions \eqref{eq:match cond 1} - \eqref{eq:match cond 5} for $\mathcal{M}''$, where $\mathcal{M}'' \lessdot_{j-1} \mathcal{M}$ is of type IA1.  This completes the proof.
\end{proof}

\brem
Note that it follows from the proof that only covering relations of types IA, IB and II are needed to generate the entire $\mathcal{W}$-set of any involution in $\mathscr{I}_n$.
\erem

\subsection{Fixed Point Free Involutions}\label{sec: fpf inv W-set}

In this subsection, $n = 2k$ is a positive even integer.

\bdfn
Let $\mathscr{I}'_n = \{ \pi' \in \mathscr{I}_n : \pi'(i) \neq i \text{ for all } 1 \leq i \leq n \}$ be the set of fixed-point free involutions of $S_n$.  It is a $M(S_n)$-set through the action defined inductively by
\begin{equation}\label{eqn:RS mon action on I'_n}
m(s_i) \cdot \pi' =
\begin{cases}
s_i \pi' s_i^{-1} & \text{if $\ell(s_i \pi' s_i^{-1}) = \ell(\pi') + 2$,} \\
\pi' & \text{otherwise.}
\end{cases}
\end{equation}
\edfn

\brem
This action is simply the restriction to $\mathscr{I}'_n$ of the action of $M(S_n)$ on $\mathscr{I}_n$.
\erem

The \emph{weak order} on $\mathscr{I}'_n$ is the restriction of the weak order on $\mathscr{I}_n$ to $\mathscr{I}'_n$.
The Hasse diagram for weak order for $\mathscr{I}'_6$ is shown in Figure \ref{fig:weak order I'_6}.  The poset $\mathscr{I}'_n$ has a bottom element $\alpha'_n = (1,2)(3,4)\cdots(n-1,n)$ and a top element $\beta'_n = w_0$.  As before, we drop the subscript $n$ when it is clear from context.

\begin{figure}[htp]
\begin{center}

\begin{tikzpicture}[scale=.4]

\node at (0,0) (a) {$(12)(34)(56)$};

\node at (-5,5) (b1) {$(13)(24)(56)$};
\node at (5,5) (b2) {$(12)(35)(46)$};

\node at (-10,10) (c1) {$(14)(23)(56)$};
\node at (0,10) (c2) {$(13)(25)(46)$};
\node at (10,10) (c3) {$(12)(36)(45)$};

\node at (-10,15) (d1) {$(15)(23)(46)$};
\node at (0,15) (d2) {$(14)(25)(36)$};
\node at (10,15) (d3) {$(13)(26)(45)$};

\node at (-10,20) (e1) {$(15)(24)(36)$};
\node at (10,20) (e2) {$(14)(26)(35)$};
\node at (0,20) (e3) {$(16)(23)(45)$};

\node at (-5,25) (f1) {$(15)(26)(34)$};
\node at (5,25) (f2) {$(16)(24)(35)$};

\node at (0,30) (g) {$(16)(25)(34)$};

\node at (-3,2.5) {$2$};
\node at (3,2.5) {$4$};

\node at (-8.6,7.5) {$1,3$};
\node at (-3.2,7.5) {$4$};
\node at (3.2,7.5) {$2$};
\node at (8.6,7.5) {$3,5$};

\node at (-10.5,12.5) {$4$};
\node at (-6,12.5) {$1$};
\node at (-0.5,12.5) {$3$};
\node at (6,12.5) {$5$};
\node at (10.5,12.5) {$2$};

\node at (-10.5,17.5) {$3$};
\node at (-8,16.5) {$5$};
\node at (-1.5,16.5) {$1,4$};
\node at (1.5,16.5) {$2,5$};
\node at (8,16.5) {$1$};
\node at (10.5,17.5) {$3$};

\node at (-8.1,22.5) {$2$};
\node at (-6,20.8) {$5$};
\node at (6,20.8) {$4$};
\node at (1,21.5) {$3$};
\node at (8.1,22.5) {$1$};

\node at (-3.6,27.5) {$1,5$};
\node at (3.6,27.5) {$2,4$};

\draw[-, thick] (a) to (b1);
\draw[-, thick] (a) to (b2);

\draw[-, thick] (b1) to (c1);
\draw[-, thick] (b1) to (c2);
\draw[-, thick] (b2) to (c2);
\draw[-, thick] (b2) to (c3);

\draw[-, thick] (c1) to (d1);
\draw[-, thick] (c2) to (d1);
\draw[-, thick] (c2) to (d2);
\draw[-, thick] (c2) to (d3);
\draw[-, thick] (c3) to (d3);

\draw[-, thick] (d1) to (e1);
\draw[-, thick] (d1) to (e3);
\draw[-, thick] (d2) to (e1);
\draw[-, thick] (d2) to (e2);
\draw[-, thick] (d3) to (e2);
\draw[-, thick] (d3) to (e3);

\draw[-, thick] (e1) to (f1);
\draw[-, thick] (e1) to (f2);
\draw[-, thick] (e2) to (f1);
\draw[-, thick] (e2) to (f2);
\draw[-, thick] (e3) to (f2);

\draw[-, thick] (f1) to (g);
\draw[-, thick] (f2) to (g);

\end{tikzpicture}

\caption{Weak order on $\mathscr{I}'_6$.}\label{fig:weak order I'_6}

\end{center}
\end{figure}

In terms of matchings, $\mathscr{I}'_n$ corresponds to matchings on $n$ vertices with no isolated vertices.  Our next
definition allows us to parameterize the maximal chains of an interval $[\alpha', \pi']$ in $\mathscr{I}'_n$.  For
$\pi' = (a_1, b_1) \cdots (a_k, b_k) \in \mathscr{I}'_n$, written in standard form, define $L'(\pi')$ to be the
number of inversions of the word $(a_1, b_1, \dots, a_k, b_k)$.  It is easily checked that $L'(\pi') = L(\pi') - k$.
Weak order on $\mathscr{I}'_n$ is a ranked poset, with rank function $L'$ \cite{RS90}.

\bdfn
The $\mathcal{W}$-set of $\pi' \in \mathscr{I}'_n$ is
$$
\mathcal{W}'(\pi') := \{ w \in S_n : m(w) \cdot \alpha' = \pi' \text{ and } \ell(w) = L'(\pi') \}.
$$
\edfn

\brem\label{rem: W-set notation}
Note that we use the notation $\mathcal{W}'(\pi)$ to denote the $\mathcal{W}$-set of an element $\pi \in
\mathscr{I}'_n$.  The reason for this is that in the proof of Corollary \ref{thm:W-set for I'_n}, we will have need
to refer to the $\mathcal{W}$-set of $\pi$ viewed as an element of the weak order poset of all involutions, for which
we reserve the notation $\mathcal{W}(\pi)$.
\erem

The computation of the $\mathcal{W}$-sets for elements of $\mathscr{I}'_n$ follows from Theorem \ref{thm:W-set for I_n}, as we now prove.

\bcor\label{thm:W-set for I'_n}
Let $n = 2k$ and let $\pi' = (a_1, b_1)(a_2, b_2) \cdots (a_k, b_k) \in \mathscr{I}'_n$ be written in standard form.  Then $\mathcal{W}'(\pi')$ consists of all $w = [w(1), w(2), \dots, w(n)]$ such that
\benum
\item for each $i$, $a_i$ occurs before $b_i$ in $w$ and no value occurs between $b_i$ and $a_i$ in $w$;\label{eq:fpf cond 1}
\item if $i < j$ and $b_i < b_j$, then $b_i$ occurs before $a_j$ in $w$.\label{eq:fpf cond 2}
\eenum
\ecor

\brem
Note that condition \eqref{eq:fpf cond 1} forces each pair of values $a_i$ and $b_i$ adjacent to each other in the one-line notation of $w$, with $a_i$ coming before $b_i$.  So one only has to consider the $k!$ possible orderings of these two-element blocks.  To produce all $W$-set elements, one can start with the permutation $w=[a_1,b_1,\hdots,a_k,b_k]$.  According to Corollary \ref{thm:W-set for I'_n} and the definition of standard form, it is always the case that $w \in \mathcal{W}'(\pi')$.  To enumerate the remaining $W$-set elements starting from $w$, one simply permutes the two-element blocks of $w$ by successively interchanging adjacent two-element blocks, but with a restriction:  \eqref{eq:fpf cond 2} says that two-element blocks $a_i,b_i$ and $a_j,b_j$ can only be interchanged if $b_i > b_j$.
\erem

\bexam
In $S_6$, $S_8$, and $S_8$, respectively, we have
\begin{align*}
\mathcal{W}'((1,6)(2,5)(3,4)) & = \{ 162534, 163425, 251634, 253416, 341625, 342516 \}, \\
\mathcal{W}'((1,6)(2,3)(4,8)(5,7)) & = \{ 16234857, 23164857, 16235748, 23165748 \}, \\
\mathcal{W}'((1,5)(2,7)(3,8)(4,6)) & = \{ 15273846, 15274638, 15462738 \}.
\end{align*}
\eexam

\begin{proof}

Consider $\alpha' = (1,2)(3,4)\cdots(n-1,n)$ as an element of $\mathscr{I}_{n}$.  Let $w^* = [2,1,4,3,\dots,n,n-1]$.  Then it follows from Theorem \ref{thm:W-set for I_n} that the $\mathcal{W}$-set of $\alpha'$ as an element of $\mathscr{I}_n$ is $\mathcal{W}(\alpha') = \{ w^* \}$.  Assume that $w'$ satisfies \eqref{eq:fpf cond 1} and \eqref{eq:fpf cond 2}; we will prove that $w' \in \mathcal{W}'(\pi')$.  Since $w'$ satisfies \eqref{eq:fpf cond 1}, it follows that $w'(i) < w'(i+1)$ for all odd $i$, so $w = w' w^*$ has length $\ell(w) = \ell(w') + \ell(w^*)$.  Moreover, conditions \eqref{eq:fpf cond 1} and \eqref{eq:fpf cond 2} of Corollary \ref{thm:W-set for I'_n} for $w'$ imply conditions \eqref{eq:match cond 1} and \eqref{eq:match cond 2} of Theorem \ref{thm:W-set for I_n} for $w$.  Since conditions \eqref{eq:match cond 3}, \eqref{eq:match cond 4}, and \eqref{eq:match cond 5} of Theorem \ref{thm:W-set for I_n} hold vacuously (as $\pi'$ has no fixed points), it follows that $w \in \mathcal{W}(\pi')$, and consequently, $w' \in \mathcal{W}'(\pi')$.

Conversely, suppose $w' \in \mathcal{W}'(\pi')$. There are no covering relations of type IA or type II in the interval $[\alpha', \pi']$ because elements of $\mathscr{I}'_n$ do not have fixed points, so their matchings do not have isolated vertices.  But in the type IB and IC covering relations, the strands containing vertex $i$ and $i+1$ are switched, and it easily checked that if $w'$ satisfies \eqref{eq:fpf cond 1} for $\mathcal{M}'$ and $\mathcal{M}' \lessdot_{i} \mathcal{M}$ is a covering relation of type IB or IC, then $s_i w'$ satisfies \eqref{eq:fpf cond 1} as well.  Since clearly the unique element $\text{id}$ of $\mathcal{W}'(\alpha')$ satisfies \eqref{eq:fpf cond 1}, it follows by induction on $L'(\pi')$ that $w' \in \mathcal{W}'(\pi')$ satisfies \eqref{eq:fpf cond 1}.  Given that every $w' \in \mathcal{W}'(\pi')$ satisfies \eqref{eq:fpf cond 1}, then $w'$ satisfies \eqref{eq:fpf cond 2} if and only if $w = w' w^*$ satisfies \eqref{eq:match cond 2} of Theorem \ref{thm:W-set for I_n}, which it does by virtue of belonging to $\mathcal{W}(\pi')$.
\end{proof}

\subsection{Involutions with Signed Fixed Points}\label{sec: sgn inv W-set}

Let $n$ be a positive integer and fix $p, q > 0$ with $p + q = n$.  Let $\mathscr{I}^{\pm}_{p,q}$ be the set of all involutions $\pi \in \mathscr{I}_n$ with an assignment of $+$ and $-$ signs to the fixed points of $\pi$ such that there are $p - q$ more $+$'s than $-$'s if $p \geq q$, or $q - p$ more $-$'s than $+$'s if $p \leq q$.  We call the quantity $p - q$ the total charge of $\pi$.  (In the literature, these combinatorial gadgets have often been referred to as ``clans".  The nomenclature can be traced back to \cite{MO90}.  They are also examples of ``decorated permutations,'' which occur in \cite{Postnikov06}.)  We denote such objects $\pi \in \mathscr{I}^{\pm}_{p,q}$ by writing the involution in cycle notation, but including the one-cycles so that we may adorn each with a $+$ or $-$ sign.  For example, $\pi = (1,6)(2,3)(4^+)(5^-)(7^+)$ is an element of $\mathscr{I}^{\pm}_{4,3}$.  Note that $p$ is equal to the number of fixed points in $\pi$ with a $+$ sign attached plus the number of two-cycles in $\pi$, while $q$ is equal to the number of fixed points in $\pi$ with a $-$ sign attached plus the number of two-cycles in $\pi$.  Without loss of generality, we will assume that $p \geq q$ in the sequel.

Just as matchings give a useful model for involutions, we use \emph{signed matchings} as a combinatorial model for elements of $\mathscr{I}^{\pm}_{p,q}$.  A signed matching is a matching $\mathcal{M}$ with an assignment of a $+$ or $-$ sign to each isolated vertex of $\mathcal{M}$.  To each $\pi \in \mathscr{I}^{\pm}_{p,q}$, we naturally assign a signed matching $\mathcal{M}_{\pi}$ that has $p - q$ more $+$'s than $-$'s.

\begin{figure}[htp]
\begin{center}

\begin{tikzpicture}[scale=.75, transform shape]

\tikzstyle{dot} = [circle, minimum width = 4 pt, fill]

\node (9) {};
\node[dot] (10) [right=of 9, label=below:$i$, label=above:$\pm$] {};
\node (11) [right=of 10] {};
\node[dot] (12) [right=of 11, label=below:$i+1$] {};
\node (13) [right=of 12] {};

\node (u9) [above=of 9] {};

\draw [-] (9.east) -- (10.west);
\draw [-] (10.east) -- (12.west);
\draw [-] (12.east) -- (13.west);

\draw[-] [in=45, out=90] (12.north) to (u9.east);

\node (6) [right=of 13] {};
\node (7) [right=of 6] {};
\node (8) [right=of 7] {};

\draw[->] (6.east) -- (8.west) node [draw=none,midway,above=1mm] {Type IA1};

\node (1) [right=of 8]{};
\node[dot] (2) [right=of 1, label=below:$i$] {};
\node (3) [right=of 2] {};
\node[dot] (4) [right=of 3, label=below:$i+1$, label=above:$\pm$] {};
\node (5) [right=of 4] {};

\node (u1) [above=of 1] {};

\draw [-] (1.east) -- (2.west);
\draw [-] (2.east) -- (4.west);
\draw [-] (4.east) -- (5.west);

\draw[-] [in=0, out=90] (2.north) to (u1.east);

\end{tikzpicture}

\begin{tikzpicture}[scale=.75, transform shape]

\tikzstyle{dot} = [circle, minimum width = 4 pt, fill]

\node (9) {};
\node[dot] (10) [right=of 9, label=below:$i$] {};
\node (11) [right=of 10] {};
\node[dot] (12) [right=of 11, label=below:$i+1$, label=above:$\pm$] {};
\node (13) [right=of 12] {};

\node (u13) [above=of 13] {};

\draw [-] (9.east) -- (10.west);
\draw [-] (10.east) -- (12.west);
\draw [-] (12.east) -- (13.west);

\draw[-] [in=135, out=90] (10.north) to (u13.west);

\node (6) [right=of 13] {};
\node (7) [right=of 6] {};
\node (8) [right=of 7] {};

\draw[->] (6.east) -- (8.west) node [draw=none,midway,above=1mm] {Type IA2};

\node (1) [right=of 8] {};
\node[dot] (2) [right=of 1, label=below:$i$, label=above:$\pm$] {};
\node (3) [right=of 2] {};
\node[dot] (4) [right=of 3, label=below:$i+1$] {};
\node (5) [right=of 4] {};

\node (u5) [above=of 5] {};

\draw [-] (1.east) -- (2.west);
\draw [-] (2.east) -- (4.west);
\draw [-] (4.east) -- (5.west);

\draw[-] [in=180, out=90] (4.north) to (u5.west);

\end{tikzpicture}

\begin{tikzpicture}[scale=.75, transform shape]

\tikzstyle{dot} = [circle, minimum width = 4 pt, fill]

\node (9) {};
\node[dot] (10) [right=of 9, label=below:$i$] {};
\node (11) [right=of 10] {};
\node[dot] (12) [right=of 11, label=below:$i+1$] {};
\node (13) [right=of 12] {};

\node (u9) [above=of 9] {};
\node (u13) [above=of 13] {};

\draw [-] (9.east) -- (10.west);
\draw [-] (10.east) -- (12.west);
\draw [-] (12.east) -- (13.west);

\draw[-] [in=135, out=90] (10.north) to (u13.west);
\draw[-] [in=45, out=90] (12.north) to (u9.east);

\node (6) [right=of 13] {};
\node (7) [right=of 6] {};
\node (8) [right=of 7] {};

\draw[->] (6.east) -- (8.west) node [draw=none,midway,above=1mm] {Type IB};

\node (1) [right=of 8] {};
\node[dot] (2) [right=of 1, label=below:$i$] {};
\node (3) [right=of 2] {};
\node[dot] (4) [right=of 3, label=below:$i+1$] {};
\node (5) [right=of 4] {};

\node (u1) [above=of 1] {};
\node (u5) [above=of 5] {};

\draw [-] (1.east) -- (2.west);
\draw [-] (2.east) -- (4.west);
\draw [-] (4.east) -- (5.west);

\draw[-] [in=0, out=90] (2.north) to (u1.east);
\draw[-] [in=180, out=90] (4.north) to (u5.west);

\end{tikzpicture}

\begin{tikzpicture}[scale=.75, transform shape]

\tikzstyle{dot} = [circle, minimum width = 4 pt, fill]

\node (9) {};
\node[dot] (10) [right=of 9, label=below:$i$] {};
\node (11) [right=of 10] {};
\node[dot] (12) [right=of 11, label=below:$i+1$] {};
\node (13) [right=of 12] {};

\node (u13) [above=of 13] {};
\node (uu13) [above=of u13] {};

\draw [-] (9.east) -- (10.west);
\draw [-] (10.east) -- (12.west);
\draw [-] (12.east) -- (13.west);

\draw[-] [in=165, out=90] (10.north) to (uu13.west);
\draw[-] [in=150, out=90] (12.north) to (u13.west);

\node (6) [right=of 13] {};
\node (7) [right=of 6] {};
\node (8) [right=of 7] {};

\draw[->] (6.east) -- (8.west) node [draw=none,midway,above=1mm] {Type IC1};

\node (1) [right=of 8] {};
\node[dot] (2) [right=of 1, label=below:$i$] {};
\node (3) [right=of 2] {};
\node[dot] (4) [right=of 3, label=below:$i+1$] {};
\node (5) [right=of 4] {};

\node (u5) [above=of 5] {};
\node (uu5) [above=of u5] {};

\draw [-] (1.east) -- (2.west);
\draw [-] (2.east) -- (4.west);
\draw [-] (4.east) -- (5.west);

\draw[-] [in=135, out=90] (2.north) to (u5.west);
\draw[-] [in=180, out=90] (4.north) to (uu5.west);

\end{tikzpicture}

\begin{tikzpicture}[scale=.75, transform shape]

\tikzstyle{dot} = [circle, minimum width = 4 pt, fill]

\node (9) {};
\node[dot] (10) [right=of 9, label=below:$i$] {};
\node (11) [right=of 10] {};
\node[dot] (12) [right=of 11, label=below:$i+1$] {};
\node (13) [right=of 12] {};

\node (u9) [above=of 9] {};
\node (uu9) [above=of u9] {};

\draw [-] (9.east) -- (10.west);
\draw [-] (10.east) -- (12.west);
\draw [-] (12.east) -- (13.west);

\draw[-] [in=30, out=90] (10.north) to (u9.east);
\draw[-] [in=15, out=90] (12.north) to (uu9.east);

\node (6) [right=of 13] {};
\node (7) [right=of 6] {};
\node (8) [right=of 7] {};

\draw[->] (6.east) -- (8.west) node [draw=none,midway,above=1mm] {Type IC2};

\node (1) [right=of 8] {};
\node[dot] (2) [right=of 1, label=below:$i$] {};
\node (3) [right=of 2] {};
\node[dot] (4) [right=of 3, label=below:$i+1$] {};
\node (5) [right=of 4] {};

\node (u1) [above=of 1] {};
\node (uu1) [above=of u1] {};

\draw [-] (1.east) -- (2.west);
\draw [-] (2.east) -- (4.west);
\draw [-] (4.east) -- (5.west);

\draw[-] [in=0, out=90] (2.north) to (uu1.east);
\draw[-] [in=45, out=90] (4.north) to (u1.east);

\end{tikzpicture}

\begin{tikzpicture}[scale=.75, transform shape]

\tikzstyle{dot} = [circle, minimum width = 4 pt, fill]

\node (9) {};
\node[dot] (10) [right=of 9, label=below:$i$] {};
\node (11) [right=of 10] {};
\node[dot] (12) [right=of 11, label=below:$i+1$] {};
\node (13) [right=of 12] {};

\draw [-] (9.east) -- (10.west);
\draw [-] (10.east) -- (12.west);
\draw [-] (12.east) -- (13.west);

\draw[-] [in=90, out=90] (10.north) to (12.north);

\node (6) [right=of 13] {};
\node (7) [right=of 6] {};
\node (8) [right=of 7] {};

\draw[->] (6.east) -- (8.west) node [draw=none,midway,above=1mm] {Type II};

\node (1) [right=of 8] {};
\node[dot] (2) [right=of 1, label=below:$i$, label=above:$\pm$] {};
\node (3) [right=of 2] {};
\node[dot] (4) [right=of 3, label=below:$i+1$, label=above:$\mp$] {};
\node (5) [right=of 4] {};

\draw [-] (1.east) -- (2.west);
\draw [-] (2.east) -- (4.west);
\draw [-] (4.east) -- (5.west);

\end{tikzpicture}

\caption{The types of covers $\mathcal{M}' \lessdot_{i} \mathcal{M}$ for signed matchings.}\label{fig:signed matching types}

\end{center}
\end{figure}
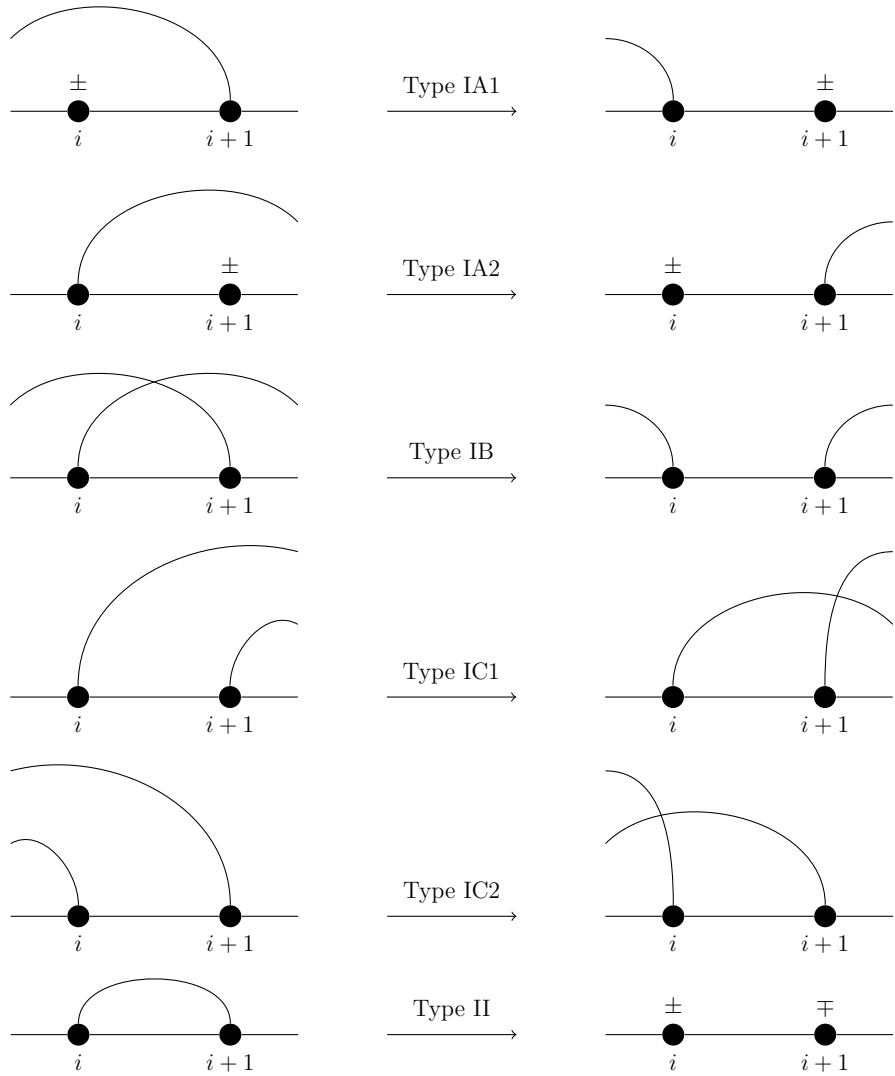

The \emph{weak order} on $\mathscr{I}^{\pm}_{p,q}$  does not correspond to a well-defined $M(S_n)$ action on $\mathscr{I}^{\pm}_{p,q}$.  (It should be noted that there is an $M(S_n)$-action on the opposite of $\mathscr{I}^{\pm}_{p,q}$ for geometric reasons, but we have chosen not to take this perspective, so as to keep our presentation entirely combinatorial.)  The covering relations for the weak order on $\mathscr{I}^{\pm}_{p,q}$, which are illustrated in Figure \ref{fig:signed matching types} in terms of signed matchings, are as follows:
\begin{itemize}
 \item Switch the endpoint of a strand with an adjacent sign so as to shorten the strand (Types IA1 and IA2);
 \item Undo a crossing occurring at consecutive vertices so as to create two disjoint strands (Type IB);
 \item Create a crossing from a nested pair of strands by crossing the ends of the strands at consecutive vertices (Types IC1 and IC2);
 \item Replace a strand of length $1$ by a pair of opposite signs (Type II).
\end{itemize}
Note that these covering relations are very similar to the covering relations for the \emph{opposite} of the weak order poset of involutions (see Figure \ref{fig:matching types}), with the only differences being the inclusion of signs for fixed points (isolated vertices).

Weak order on $\mathscr{I}^{\pm}_{p,q}$ is a graded poset with length function $L^{\pm}(\pi) = pq - L(\pi)$, where by $L(\pi)$ we refer to the length of the underlying ordinary involution, as defined in Section \ref{sec: inv W-set} \cite{RS90}.  Note that $\mathscr{I}^{\pm}_{p,q}$ has a unique minimal element $\alpha := (1,n)(2,n-1)\cdots(q,n+1-q)((q+1)^+)\cdots((n-q)^+)$, but has many maximal elements.  Indeed, the maximal elements correspond to the identity permutation with a choice of $p$ values to be assigned $+$'s and $q$ values to be assigned $-$'s.  Thus, there are ${n \choose p} = {n \choose q}$ maximal elements.  Weak order for $\mathscr{I}^{\pm}_{2,2}$ is illustrated in Figure \ref{fig:weak order I^+-_2,2}.

\begin{figure}[htp]
\begin{center}

\begin{tikzpicture}[xscale=.23, yscale=.5, every node/.style={scale=.85}]

\node at (0,0) (a) {$(14)(23)$};

\node at (-10,5) (b1) {$(14)(2^+)(3^-)$};
\node at (0,5) (b2) {$(13)(24)$};
\node at (10,5) (b3) {$(14)(2^-)(3^+)$};

\node at (-20,10) (c1) {$(1^+)(24)(3^-)$};
\node at (-10,10) (c2) {$(13)(2^-)(4^+)$};
\node at (0,10) (c3) {$(12)(34)$};
\node at (10,10) (c4) {$(13)(2^-)(4^+)$};
\node at (20,10) (c5) {$(1^-)(24)(3^+)$};

\node at (-25,15) (d1) {$(1^+)(23)(4^-)$};
\node at (-15,15) (d2) {$(1^+)(2^-)(34)$};
\node at (-5,15) (d3) {$(12)(3^+)(4^-)$};
\node at (5,15) (d4) {$(12)(3^-)(4^+)$};
\node at (15,15) (d5) {$(1^-)(2^+)(34)$};
\node at (25,15) (d6) {$(1^-)(23)(4^+)$};

\node at (-25,20) (e1) {$1^+ 2^+ 3^- 4^-$};
\node at (-15,20) (e2) {$1^+ 2^- 3^+ 4^-$};
\node at (-5,20) (e3) {$1^+  2^- 3^- 4^+$};
\node at (5,20) (e4) {$1^- 2^+ 3^+ 4^-$};
\node at (15,20) (e5) {$1^- 2^+ 3^- 4^+$};
\node at (25,20) (e6) {$1^- 2^- 3^+ 4^+$};

\node at (-6.5,2.5) {$2$};
\node at (1.25,2.5) {$1,3$};
\node at (6.5,2.5) {$2$};

\node at (-16.5,7.5) {$1$};
\node at (-11,7.5) {$3$};
\node at (-1,7.5) {$2$};
\node at (11,7.5) {$3$};
\node at (16.5,7.5) {$1$};

\node at (-22,10.8) {$3$};
\node at (-18.5,10.8) {$2$};
\node at (-12,11.2) {$1$};
\node at (-8.5,10.8) {$2$};
\node at (-4,10.8) {$1$};
\node at (-1,12) {$3$};
\node at (1,12) {$3$};
\node at (4,10.8) {$1$};
\node at (8.5,10.8) {$2$};
\node at (12,11.2) {$1$};
\node at (18.5,10.8) {$2$};
\node at (22,10.8) {$3$};

\node at (-26, 15.8) {$2$};
\node at (-22, 15.8) {$2$};
\node at (-16, 15.8) {$3$};
\node at (-12, 15.8) {$3$};
\node at (-8, 15.8) {$1$};
\node at (-2, 15.8) {$1$};
\node at (2, 15.8) {$1$};
\node at (8, 15.8) {$1$};
\node at (12, 15.8) {$3$};
\node at (16, 15.8) {$3$};
\node at (22, 15.8) {$2$};
\node at (26, 15.8) {$2$};

\draw[-, thick] (a) to (b1);
\draw[-, thick] (a) to (b2);
\draw[-, thick] (a) to (b3);

\draw[-, thick] (b1) to (c1);
\draw[-, thick] (b1) to (c2);
\draw[-, thick] (b2) to (c3);
\draw[-, thick] (b3) to (c4);
\draw[-, thick] (b3) to (c5);

\draw[-, thick] (c1) to (d1);
\draw[-, thick] (c1) to (d2);
\draw[-, thick] (c2) to (d1);
\draw[-, thick] (c2) to (d3);
\draw[-, thick] (c3) to (d2);
\draw[-, thick] (c3) to (d3);
\draw[-, thick] (c3) to (d4);
\draw[-, thick] (c3) to (d5);
\draw[-, thick] (c4) to (d4);
\draw[-, thick] (c4) to (d6);
\draw[-, thick] (c5) to (d5);
\draw[-, thick] (c5) to (d6);

\draw[-, thick] (d1) to (e1);
\draw[-, thick] (d1) to (e2);
\draw[-, thick] (d2) to (e2);
\draw[-, thick] (d2) to (e3);
\draw[-, thick] (d3) to (e2);
\draw[-, thick] (d3) to (e4);
\draw[-, thick] (d4) to (e3);
\draw[-, thick] (d4) to (e5);
\draw[-, thick] (d5) to (e4);
\draw[-, thick] (d5) to (e5);
\draw[-, thick] (d6) to (e5);
\draw[-, thick] (d6) to (e6);

\end{tikzpicture}

\caption{Weak order on $\mathscr{I}^{\pm}_{2,2}$.}\label{fig:weak order I^+-_2,2}

\end{center}
\end{figure}

Since $\mathscr{I}^{\pm}_{p,q}$ is not equipped with an $M(S_n)$-action, we define the
$\mathcal{W}$-sets directly in terms of chains of signed matchings.

\bdfn
The $\mathcal{W}$-set of $\pi \in \mathscr{I}^{\pm}_{p,q}$ consists of all $w \in S_n$ of length $\ell(w) = L^{\pm}(\pi)$ such that for some (equivalently, any) reduced decomposition $s_{i_k} \cdots s_{i_2} s_{i_1}$,
$$
\alpha \lessdot_{i_1} \pi_1 \lessdot_{i_2} \pi_2 \lessdot_{i_3} \dots \lessdot_{i_{k-1}} \pi_{k-1} \lessdot_{i_k} \pi
$$
for some $\pi_1, \pi_2, \dots, \pi_{k-1} \in \mathscr{I}^{\pm}_{p,q}$.
\edfn

\begin{Algorithm}\label{alg:+-}
We now describe an algorithm which generates the $\mathcal{W}$-set of an element $\pi \in \mathscr{I}^{\pm}_{p,q}$.  Initialize $A_0 = [n]$.  We will inductively construct a permutation $w \in S_n$, which we prove in Theorem \ref{thm:W-set for sgn inv} belongs to $\mathcal{W}(\pi)$.  Let $\mathcal{M} := \mathcal{M}_{\pi}$ be the signed matching associated to $\pi$.  Assume that inductively, the values of the first and last $i-1$ values of $w$ have already been determined and let $A_{i-1} := [n] \setminus \{ w(1), \dots, w(i-1), w(n+2-i), \dots, w(n) \}$, which has cardinality $n - 2(i-1)$.  Then choose either (1) a strand $\{ a < b \}$ of $\mathcal{M}$ with $a, b \in A_{i-1}$ and such that the strand is not nested inside any strand $\{ c < d \}$ of $\mathcal{M}$ with $c, d \in A_{i-1}$; or (2) two isolated vertices $a < b$ of opposite sign that are adjacent in $A_{i-1}$ and also are not nested inside any other strand $\{ c < d \}$ of $\mathcal{M}$ with $c, d \in A_{i-1}$.  In case (1), set $w(i) = a$ and $w(n+1-i) = b$; in case (2), set $w(i) = b$ and $w(n+1-i) = a$.  In either case, set $A_i = A_{i-1} \setminus \{a,b\}$.  Continue until $A_k$ consists of only isolated vertices of $\mathcal{M}$ of the same sign.  Then define the remaining as yet undetermined middle values of $w$ to be the elements of $A_k$ in increasing order.  Let $\mathcal{A}(\pi)$ denote the set of all $w \in S_n$ obtained by this algorithm for some sequence of choices.
\end{Algorithm}

\bexam
\begin{align*}
\mathcal{A}((1^+)(2^-)(3^+)(4^+)) & = \{2341, 3142\}, \\
\mathcal{A}((1,6)(2,3)(4^+)(5^-)(7^+)) & = \{ 1257436, 1274536, 1527346, 1724356 \}, \\
\mathcal{A}((1,4)(2,6)(3,5)) & = \{ 123564, 213546, 231456 \}, \\
\mathcal{A}((1^+)(2^-)(3^-)(4^+)(5^+)(6^+)) & = \{ 245631, 425613, 451623 \}.
\end{align*}
\eexam

\bthm\label{thm:W-set for sgn inv}
For any $\pi \in \mathscr{I}^{\pm}_{p,q}$, $\mathcal{W}(\pi) = \mathcal{A}(\pi)$.
\ethm

\begin{proof}
We prove the theorem using induction on $L^{\pm}(\pi)$.  If $L^{\pm}(\pi) = 0$, then $\pi = \alpha$.  Following Algorithm \ref{alg:+-} to produce $w \in \mathcal{A}(\pi)$ requires us to set $w(1) = 1$, $w(n) = n$, then $w(2) = 2$, $w(n-1) = n-1$, and so on until we reach $w(q) = q$, $w(n+1-q) = n+1-q$.  Then the second stage of the algorithm requires us to set $w(q+1) = q+1, \dots, w(n-q) = n-q$.  Thus $w = \text{id}$ and $\mathcal{A}(\pi) = \{ \text{id} \} = \mathcal{W}(\pi)$.

Now fix $\pi$ with $L^{\pm}(\pi) > 0$, with $\mathcal{M}_{\pi}$ its associated signed matching.  We first show that $\mathcal{W}(\pi) \subseteq \mathcal{A}(\pi)$.  Let $w \in \mathcal{W}(\pi)$ be given.  By definition, there exists $\pi' \in \mathscr{I}^{\pm}_{p,q}$ (with associated signed matching $\mathcal{M}_{\pi'}$) such that $\pi' \lessdot_{j} \pi$ for some $1 \leq j < n$ and $w' = s_j w$ belongs to $\mathcal{W}(\pi')$.  In particular, $\ell(w') = \ell(w) - 1$.  By our induction hypothesis, $w' \in \mathcal{A}(\pi')$.  So $w'$ is obtained by carrying out Algorithm \ref{alg:+-} for $\mathcal{M}_{\pi'}$, say using steps $S_1',\hdots,S_k'$.  We claim that if we carry out steps $S_1,\hdots,S_k$ of Algorithm \ref{alg:+-} for $\mathcal{M}_{\pi}$ which coincide with the steps $S_1',\hdots,S_k'$ except with $j$ and $j+1$ interchanged, then $w$ is produced, thus proving that $w \in \mathcal{A}(\pi)$.  Note that $\mathcal{M}_{\pi}$ is obtained from $\mathcal{M}_{\pi'}$ by swapping the vertices $j$ and $j+1$, unless $j$ and $j+1$ were matched in $\mathcal{M}_{\pi'}$, in which case that matching is deleted in $\mathcal{M}_{\pi}$ and the resulting two isolated vertices are given opposite signs.  So if a strand or two adjacent isolated vertices of $\mathcal{M}_{\pi}$ are nested in a strand in $A_i$, then the same is true for the corresponding strand or two adjacent isolated vertices of $\mathcal{M}_{\pi'}$.  Indeed, the only difference in nesting between the two matchings occurs if the covering relation is of type IC, and in that case, $\mathcal{M}_{\pi'}$ has all of the nestings of $\mathcal{M}_{\pi}$ plus an additional one.  Thus $S_1,\hdots,S_k$ is a valid sequence of steps to apply to the matching $\mathcal{M}_{\pi}$. Moreover, for any type I covering relation, the effect of carrying out $S_1,\hdots,S_k$ for $\mathcal{M}_{\pi}$ is to create $w$ that agrees with $w'$ except that $j$ and $j+1$ are interchanged, i.e. $w = s_j w'$.  The same is true for a type II covering relation because in that case, the isolated vertices $j$ and $j+1$ of $\mathcal{M}_{\pi}$ have opposite sign, so whereas $j$ and $j+1$ are placed in order in $w'$ in one of the steps $S_i'$ for $\mathcal{M}_{\pi'}$, they are placed in reverse order in $w$ in the corresponding step $S_i$.  So again $w = s_j w'$.

We now show that $\mathcal{A}(\pi) \subseteq \mathcal{W}(\pi)$.  Let $w \in \mathcal{A}(\pi)$ be given.  Since $\pi \neq \alpha$, $w \neq \text{id}$, so in particular, $w$ has at least one left descent.  (Recall that a left descent of $w$ is a number $1 \leq j \leq n-1$ such that $j+1$ occurs before $j$ in $w$.)  We first claim that if $j$ is any left descent of $w$, then there exists $\pi' \in \mathscr{I}^{\pm}_{p,q}$ such that $\pi' \lessdot_j \pi$.  Indeed, when carrying out Algorithm \ref{alg:+-} for $\mathcal{M}$ to produce $w$, we can make $j$ a left descent in several ways.  For one, we could first encounter $j$ as the right vertex of a strand and then encounter $j+1$ as an isolated vertex.  This corresponds to a covering relation of type IA1.  Similarly, we could first encounter $j+1$ as the left vertex of a strand and then encounter $j$ as an isolated vertex, corresponding to a covering relation of type IA2.  Or both $j$ and $j+1$ could belong to distinct strands.  In order for $j$ to be a left descent in this case, we must either choose the strand with $j$ in it first and have $j$ as the right vertex of its strand (type IB or IC2), or we must choose the strand with $j+1$ in it first and have $j+1$ as the left vertex in its strand (type IB or IC1).  Finally, if $j$ and $j+1$ are isolated vertices, then they must have opposite sign and be chosen at the same time (type II).  Indeed, if $j$ is matched up with an oppositely signed isolated vertex $k < j$ and $j+1$ matched up with an oppositely signed isolated vertex $l > j+1$, then they would appear in the order $\hdots j \dots l \dots (j+1) \dots k \hdots$ or $\hdots l \dots j \dots k \dots j+1 \hdots$ in $w$, contradicting the fact that $j$ is assumed to be a left descent of $w$.  A similar argument applies if either $j$ or $j+1$ is an unmatched isolated vertex from the second stage of Algorithm \ref{alg:+-}.  Note that in all the cases, there exists $\pi' \lessdot_j \pi$.  More precisely, in type I, $\pi' = s_j \pi s_j$ (with signs of fixed points determined in the obvious manner) and in type II, $\pi'=  s_j \pi$.  That $\pi'$ is covered by $\pi$ in all cases (and not the other way around) follows from the explicit description of the covering relations.

Now, for the $\pi'$ whose existence we have established, let $w' = s_j w$.  We claim that $w' \in \mathcal{A}(\pi')$.  Then, by induction, $w' \in \mathcal{W}(\pi')$ and so $w \in \mathcal{W}(\pi)$, completing the proof.  To establish our claim, we again argue that we may carry out the same sequence of choices in Algorithm \ref{alg:+-} for $\mathcal{M}'$ as was made to produce $w$ by applying Algorithm \ref{alg:+-} for $\mathcal{M}$, except that $j$ and $j+1$ must be interchanged as appropriate.  The key observation is in type IC, where $\mathcal{M}'$ has an additional nesting that does not occur in $\mathcal{M}$.  Without loss of generality, assume we are in type IC1.  Thus $\mathcal{M}$ matches $j$ and some $k > j+1$, as well as $j+1$ and some $l > k$.  So $\mathcal{M}'$ matches $j$ and $l$, as well as $j+1$ and $k$.  Because of the nesting, the order of these four elements in $w'$ for \emph{any} sequence of choices made will be $\hdots j \hdots j+1 \hdots k \hdots l \hdots$.  In $w$, where there is no nesting between these two cycles, the four elements can \emph{a priori} appear in either the order $\hdots j+1 \hdots j \hdots k \hdots l \hdots$ or $\hdots j \hdots j+1 \hdots l \hdots k \hdots$.  But for our particular $w$, the second possibility is ruled out because we are assuming that $j$ is a left descent of $w$.  Thus, the sequence of choices in Algorithm \ref{alg:+-} for $\mathcal{M}$ to produce $w$ will result in the strand $\{ j+1 < l \}$ being chosen before the strand $\{ j < k \}$.  Consequently, we may make the same sequence of choices in Algorithm \ref{alg:+-} for $\mathcal{M}'$ with the strand $\{ j < l \}$ in $\mathcal{M}'$ being chosen in place of the strand $\{ j+1 < l \}$ in $\mathcal{M}$ and the strand $\{ j+1 < k \}$ in $\mathcal{M}'$ chosen in place of the strand $\{ j < k \}$ in $\mathcal{M}$.  This results in $w' = s_j w$, and so $w' \in \mathcal{A}(\pi')$, as claimed.
\end{proof}

\brem
For a special class of clans $\pi$ (those for which no two strands cross), $\mathcal{W}(\pi)$ has significance in certain Schubert calculus problems.  Indeed, for such a $\pi$, it is shown in \cite{Wyser-JAC} that in the Chow ring of the complete flag variety $GL_n / B$, we have the following identity among Schubert classes, for permutations $u(\pi)$ and $v(\pi)$ explicitly defined in terms of the combintorics of $\pi$:
\[ [X^{u(\pi)}] \cdot [X^{v(\pi)}] = \displaystyle\sum_{w \in \mathcal{W}(\pi)} [X^{w^{-1}}]. \]
Thus Theorem \ref{thm:W-set for sgn inv} enhances this result by specifying which $w$ occur on the right hand side more explicitly.
\erem

\vspace{.5cm}

\noindent \textbf{Acknowledgement.}
The first and the second authors are partially supported by N.S.A. Grant H98230-14-1-0142.
The first author is partially supported by the Louisiana Board of Regents Research and Development Grant 549941C1.
The third author is supported by NSF International Research Fellowship 1159045, and is hosted by Institut Fourier in
Grenoble.

\bibliography{References}
\bibliographystyle{plain}

\end{document}